\documentclass[12pt,twoside]{amsart}
\usepackage{euler}
\usepackage{epsfig}
\usepackage{latexsym}
\usepackage{graphicx}
 \def\Dj{\hbox{D\kern-.73em\raise.30ex\hbox{-} \raise-.30ex\hbox{}}}
 \def\dj{\hbox{d\kern-.33em\raise.80ex\hbox{-} \raise-.80ex\hbox{\kern-.40em}}}

\begin{document}

\textit{  $ \quad \;$      Atti Semin. Mat. Fis. Univ. Modena Reggio
Emilia}, 56 (2008-2009),  95-111.

\title[RIEMANN'S FUNCTIONAL EQUATION ALTERNATIVE FORM]{$\mathbb{\quad\;AN\;\;ALTERNATIVE\;\;FORM\;\;OF\;\;THE\;\;FUNCTIONAL\;\;EQUATION}\newline\;\;
\mathbb{FOR\;\;RIEMANN'S\;\;ZETA\;\;FUNCTION}$}

\vspace{8mm}

\author{$\mathrm{ANDREA\;\;OSSICINI}$}

\maketitle

\vspace{10mm}

\baselineskip=0.20in

\vspace{10mm}

\baselineskip=0.17in \noindent {\small {\bf Abstract.} In this paper
we present a simple method for deriving an alternative form of the
functional equation for Riemann's Zeta function. The connections
between some functional equations obtained implicitly by Leonhard
Euler in his work \textit{"Remarques sur un beau rapport entre les
series des puissances tant directes que reciproques"} in Memoires de
l'Academie des Sciences de Berlin 17, (1768), permit to define a
special function, named $A\left( s \right)$, which is fully
symmetric and is similar to Riemann's $\xi $ function\footnote{$\xi
\left( s \right) = \prod\limits {\left( {s \mathord{\left/
{\vphantom {s 2}} \right. \kern-\nulldelimiterspace} 2} \right)}
\;\left( {s - 1} \right)\;\pi ^{{ - s} \mathord{\left/ {\vphantom {{
- s} 2}} \right. \kern-\nulldelimiterspace} 2}\;\zeta \left( s
\right)$ }. To be complete we find several integral representations
of the $A\left( s \right)$ function  and as a direct consequence of
the second integral representation we obtain also an analytic
continuation of the same function using an identity of Ramanujan.

\vspace{2mm}

2000\textit{ Mathematics Subject Classification: {Primary 11M35;
Secondary 11B68, 11M06.}}

\vspace{5mm} Keywords: Riemann Zeta, Dirichlet Beta, Riemann
Hypothesis, series representations.

\baselineskip=0.20in

\vspace{10mm}

\begin{center}
1. INTRODUCTION
\end{center}

\vspace{6mm}

Formulae (1) and (2) below can be found in the chapter devoted to
Euler's Gamma function in [4].

\vspace{3mm}

These are namely two functional equations for the Eulerian Zeta and
for the alternating Zeta, connected with the odd numbers, best known
as Dirichlet's Beta function and Catalan's Beta function, see [4,
pag. 35, formulae (24) and (29)].

\vspace{3mm}

Both of them were discovered, over 100 years before G.F.B. Riemann
and O. Schl\"{o}milch [8, notes on chapter II], by L. Euler in 1749
and published in 1768 in "Memoires de l'Academie des Sciences de
Berlin 17", with the title $  $ of \textit{ ``Remarques sur un beau
rapport entre les series des puissances tant directes que
reciproques''}\footnote{ Leonhardi Euleri\textit{, Opera Omnia}:
Series 1, Volume 15, pp. 70 - 90 \par }.

\newpage

The former gives, actually, an analytic extension to the complex
half-plane with $\Re \left( s \right) <  1$:

\begin{equation}
\label{eq1} \zeta \left( {1 - s} \right) = 2\left( {2\pi } \right)^{
- s}\Gamma \left( s \right)  \cos \left( {\frac{\pi s}{2}}
\right)\zeta \left( s \right).
\end{equation}

\vspace{3mm}

Here, $\Gamma$ denotes Euler's Gamma function.

\vspace{3mm}

Let us remember that the Riemann Zeta function $\zeta \left( s
\right)$ is defined by [15, pp. 96-97,see Section 2.3]:

\[
\zeta \left( s \right): = \left\{ {\begin{array}{l}
 \quad \sum\limits_{n = 1}^\infty {\frac{1}{n^s}} = \frac{1}{1 - 2^{ -
s}}\sum\limits_{n = 1}^\infty {\frac{1}{\left( {2n - 1}
\right)^s}\quad
\quad \left( {\Re \left( s \right) > 1} \right)} \\
 \\
 \;\quad \left( {1 - 2^{1 - s}} \right)^{ - 1}\;\sum\limits_{n = 1}^\infty
{\frac{\left( { - 1} \right)^{n - 1}}{n^s}\quad \quad \;\left( {\Re
\left( s
\right) > 0\,;\,s \ne 1} \right)} \\
 \end{array}} \right.
\]

\vspace{3mm}

\noindent which can indeed be analytically continued to the whole
complex $s$ plane except for a simple pole at $s = 1$ with residue
1.

\vspace{3mm}

The Riemann Zeta function $\zeta \left( s \right)$ plays a central
role in the applications of complex analysis to number theory.

\vspace{3mm}

The number-theoretic properties of $\zeta \left( s \right)$ are
exhibited by the following result as \textit{Euler's product
formula}, which gives a relationship between the set of primes and
the set of positive integers:

\[
\zeta \left( s \right) = \prod\limits_p {\left( {1 - p^{ - s}}
\right)} ^{ - 1}\quad \left( {\Re \left( s \right) > 1} \right),
\]

\noindent where the product is taken over all primes.

\vspace{3mm} It is an analytic version of the fundamental theorem of
arithmetic, which states that every integer can be factored into
primes in an essentially unique way.

\vspace{2mm}

Euler used this product to prove that the sum of the reciprocals of
the primes diverges.

\vspace{3mm}

The latter, on the contrary, gives an analytic extension on the
whole complex plane of the $L\left( s \right)$ function, that is
Dirichlet's $L$ function for the nontrivial character modulo 4,
which was later denoted by $L\left( {s,\chi _4 } \right)$ by other
authors:

\vspace{3mm}

\begin{equation}
\label{eq2} L\left( {1 - s} \right) = \left( {\frac{2}{\pi }}
\right)^s \Gamma \left( s \right) \sin \left( {\frac{\pi s}{2}}
\right)L\left( s \right).
\end{equation}

\newpage

The $L\left( s \right)$ function was defined and used by Euler,
practically for $\Re \left( s \right) > 0$ with the expression:

\vspace{3mm}

\[
L\left( s \right) = \sum\limits_{n = 0}^\infty {\;\frac{\left( { -
1} \right)^n}{\left( {2n + 1} \right)^s}}.
\]

\vspace{3mm}

It does not possess any singular point.

The $L\left( s \right)$  function is also connected to the theory of
primes [5] which may perhaps be best summarized by

\vspace{1mm}

\[
L\left( s \right) = \prod\limits_{p{\kern 1pt} \equiv {\kern 1pt}
1{\kern 1pt} \bmod {\kern 1pt} 4} {\left( {1 - p^{ - s}} \right)} ^{
- 1} \cdot \prod\limits_{p{\kern 1pt} \equiv {\kern 1pt} 3{\kern
1pt} \bmod {\kern 1pt} 4} {\left( {1 - p^{ - s}} \right)} ^{ - 1} =
\prod\limits_{p\,\,odd{\kern 1pt} } {\left( {1 - \left( { - 1}
\right)^{\frac{p - 1}{2}}p^{ - s}} \right)} ^{ - 1},
\]

\vspace{3mm}

\noindent where the products are taken over primes and the
rearrangement of factors is permitted because of an absolute
convergence.

\vspace{3mm}

Among the properties of these functions, we will limit ourselves to
report the following integral representations [4, pp. 32 and 35,
formulae (4),(5) and (28)]:

\[
\zeta \left( s \right)\,\; = \;\;\frac{1}{\Gamma \left( s
\right)}\int\limits_0^\infty {\frac{t^{s - 1}}{e^t - 1}dt} \quad
\left( {\Re \left( s \right) > 1} \right)
\]

\vspace{3mm}

\[
\zeta \left( s \right)\,\; = \;\;\frac{1}{\left( {1 - 2^{1 - s}}
\right)\,\Gamma \left( s \right)}\int\limits_0^\infty {\frac{t^{s -
1}}{e^t + 1}dt} \quad \left( {\Re \left( s \right) > 0} \right)
\]

\vspace{3mm}

\[
L\left( s \right)\,\; = \;\;\frac{1}{\,\Gamma \left( s
\right)}\int\limits_0^\infty {\frac{t^{s - 1}}{e^t + e^{ - t}}dt}
\quad \left( {\Re \left( s \right) > 0} \right).\quad
\]

\vspace{3mm}

The first integral representation, due to Abel [1], is the key
ingredient of the first proof of Bernhard Riemann [14], that is the
proof of his classical functional equation for $\zeta \left( s
\right)$.

\vspace{3mm}

Riemann obtained it by making the change of variable $x = tn$ in the
definition (via the integral) of the Gamma function, also called
Eulerian integral of second kind, and then summing for all $n \ge
1$, as shown by the following sequence of formulae, for $\Re \left(
s \right) > 1$:

\[
\Gamma \left( s \right)\,\; = \;\int\limits_0^\infty {\frac{x^{s -
1}}{e^x}dx,} \quad \,\frac{\Gamma \left( s \right)}{n^s}\; =
\;\int\limits_0^\infty {\frac{t^{s - 1}}{e^{tn}}dt,} \quad and \quad
\Gamma \left( s \right)\zeta \left( s \right)\,\; =
\;\int\limits_0^\infty {\frac{t^{s - 1}}{e^t - 1}dt} .\quad
\]

\vspace{10mm}

\begin{center}
2. THE FUNCTIONAL EQUATION FOR THE RIEMANN ZETA FUNCTION
\end{center}

\vspace{4mm}

The functional equation for the Riemann Zeta function is shown in
the asymmetrical formulation:

\[
\zeta \left( s \right) = 2^s\pi ^{s - 1}\Gamma \left( {1 - s}
\right)\sin \left( {\frac{s\pi }{2}} \right)\zeta \left( {1 - s}
\right)
\]

\noindent or in the symmetrical formulation:

\[
\Gamma \left( {\frac{s}{2}} \right)\pi ^{{ - s} \mathord{\left/
{\vphantom {{ - s} 2}} \right. \kern-\nulldelimiterspace} 2}\zeta
\left( s \right) = \Gamma \left( {\frac{1 - s}{2}} \right) \;\pi ^{
- {\left( {1 - s} \right)} \mathord{\left/ {\vphantom {{\left( {1 -
s} \right)} 2}} \right. \kern-\nulldelimiterspace} 2}\zeta \left( {1
- s} \right)
\]

\vspace{3mm}

Euler verified this relation exactly for all integer values of $s$
and numerically to great accuracy for many fractional values as
well.

\vspace{3mm}

Naturally Euler worked for integer values of $s$ with what we now
call Abel Summation (see [8] and also [18]).

\vspace{3mm}

If $\sum\limits_{n = 0}^\infty {a_n } $ is a series such that
$\;\sum\limits_{n = 0}^\infty {a_n z^n} $ converges inside the unit
disk, we shall say that $\sum\limits_{n = 0}^\infty {a_n } $ is
\textit{Abel summable} to the value $s$ if
\[
\mathop {\lim }\limits_{0 \prec x \prec 1,x \to 1^ - }
\;\sum\limits_{n = 0}^\infty {a_n x^n} = s
\]

\vspace{3mm}

In particular, if the sum $f\left( z \right)$ of the power series
extends analytically to a domain containing $z = 1$, we can take
$f\left( 1 \right)$ as value of the sum: this is what Euler did.

\vspace{3mm}

Practically Euler proved, for any integer $m \ge 2$, that

\[
\mathop {\lim }\limits_{x \to 1^ - } \frac{\eta \left( {1 - m,x}
\right)}{\eta \left( {m,x} \right)}: = \left\{ {\begin{array}{l}
 \quad - \frac{\left( {2^m - 1} \right)}{\pi ^m\left( {2^{m - 1} - 1}
\right)}\Gamma \left( m \right)\,\cos \left( \frac{\pi
\,m}{2}\right)\quad
if\;m\;is\;even, \\
 \\
 \;\quad \;0\quad \quad \quad \quad \quad \quad \quad \quad \quad \quad
\quad\quad\;if\;m\;is\;odd \\
 \end{array}} \right.
\]

\noindent were the associated power series is $\eta \left( {m,x}
\right) = \;\sum\limits_{n = 1}^\infty {\frac{\left( { - 1}
\right)^{n - 1}}{n^m}x^n.\quad \quad \;} $

\vspace{3mm}

He formally replaced $\mathop {\lim }\limits_{x \to 1^ - } \;\eta
\left( {m,x} \right)\;$ by $\;\eta \left( m \right)$,\; that is
$\sum\limits_{n = 1}^\infty {\frac{\left( { - 1} \right)^{n -
1}}{n^m} = \left( {1 - 2^{1 - m}} \right)\zeta \left( m
\right).\quad \quad \;} $

\vspace{3mm}

However he did not know how to prove this intriguing assertion for
all real numbers and in 1859 Riemann was the first to indicate that
the functional equation for $\zeta \left( s \right)$ is true and
more precisely he provided two different proofs of his classical
functional equation [14].

\newpage

Riemann is also the first to use two of the basic identities of the
Gamma function, that are Euler's complement formula and Legendre's
duplication formula (the explicit formulae are given below), to
rewrite the asymmetrical formulation in the symmetrical formulation
[14].

\vspace{3mm}

Nowadays, many proofs of this important result exist.

\vspace{3mm}

Later, mathematicians like Hardy, Siegel and others enriched the
list of proofs.

\vspace{3mm}

For example in the classical treatise of Titchmarsh [17, pp. 16-27]
seven different methods are presented .

\vspace{3mm}

Many proofs are based on Poisson's summation formula and other
proofs are unnecessarily long or conceptually difficult.

\vspace{3mm}

Recently, two authors [12] have presented a short proof of Riemann's
functional equation, based upon Poisson summation.

\vspace{3mm}

This elegant and powerful technique was used to derive a new and
simple proof of Lipschitz summation formula and as a direct
consequence of this the authors obtained an easy proof of the
functional equation for $\zeta \left( s \right)$.

\vspace{6mm}

The main point of the paper is that Hurwitz's relation [15, pp.
89-90, formulae (6) and (7)]:

\[
\zeta \left( {1 - s,a} \right) = \frac{\Gamma \left( s
\right)}{\left( {2\pi } \right)^s}\left\{ {\,e^{{ - \pi \,i\,s}
\mathord{\left/ {\vphantom {{ - \pi \,i\,s} 2}} \right.
\kern-\nulldelimiterspace} 2}F\left( {s,a} \right) + e^{{\pi \,i\,s}
\mathord{\left/ {\vphantom {{\pi \,i\,s} 2}} \right.
\kern-\nulldelimiterspace} 2}F\left( {s, - a} \right)} \right\}
\]

\vspace{3mm}

\noindent can be obtained as a conceptually simple corollary of
Lipschitz summation.

\vspace{3mm}

Notice that $F\left( {s,a} \right)$, which is often referred to as
the `periodic (or Lerch) Zeta function', is defined by the following
relation:

\[
F\left( {s,a} \right) = \sum\limits_{n = 1}^\infty {\frac{e^{2\pi
\,i{\kern 1pt} n{\kern 1pt} {\kern 1pt} a}}{n^s}}
\]

\vspace{2mm}

\noindent and that when $a = 1$ in Hurwitz's relation, $F\left(
{s,1} \right) = F\left( {s, - 1} \right) = \zeta \left( s \right)$,
so that Riemann functional equation (1) follows directly from that
by using $e^{ - \frac{\pi {\kern 1pt} i{\kern 1pt} {\kern 1pt}
s}{2}} + e^{\frac{\pi {\kern 1pt} i{\kern 1pt} {\kern 1pt} s}{2}} =
2\cos \left( {\frac{\pi {\kern 1pt} s}{2}} \right)$.

\vspace{3mm}

Here we briefly present one proof of the functional equation
extracted from [17] and due to Hardy and we will see that the
functional equation is strongly related to Fourier's series.

\vspace{3mm}

The starting point of Hardy is not the function $\zeta \left( s
\right)$, but the function:

\[
\;\sum\limits_{n = 1}^\infty {\frac{\left( { - 1} \right)^{n -
1}}{n^s} = \left( {1 - 2^{1 - s}} \right)\zeta \left( s
\right).\quad }
\]

\vspace{3mm}

This Dirichlet series is convergent for all positive values of $s$,
and so, by a general theorem on the convergence of Dirichlet series,
it is convergent for all values of $s$, with $\Re \left( s \right)
> 0$.

\vspace{3mm}

Here the pole of $\zeta \left( s \right)$ at $s = 1$ is cancelled by
the zero of the other factor.

\vspace{3mm}

Hardy's proof runs as follows. Let

\[
f\left( x \right)\; = \sum\limits_{n = 0}^\infty {\frac{\sin \left(
{2n + 1} \right)x}{2n + 1}.}
\]

\vspace{5mm}

This series is boundedly convergent and

\[
\;f\left( x \right) = \left( { - 1} \right)^m\frac{1}{4}\pi \quad
for\quad m\pi < x < \left( {m + 1} \right)\pi \quad \left( {m =
0,\;1,\;2,...} \right).
\]

\vspace{3mm}

Multiplying by $\;x^{s - 1}\; \left( {0 < s < 1} \right)$, and
integrating over$\;\left( {0,\infty } \right)$, we obtain

\[
\frac{1}{4}\pi \sum\limits_{m = 0}^\infty {\left( { - 1} \right)^m }
\,\int\limits_{m\pi }^{\left( {m + 1} \right)\pi } {x^{s - 1}dx =
\Gamma \left( s \right)\sin \left( {\frac{s\pi }{2}}
\right)\sum\limits_{n = 0}^\infty {\frac{1}{\left( {2n + 1}
\right)^{s + 1}}} = }
\]

\[ = \Gamma \left( s \right)\sin \left( {\frac{s\pi
}{2}} \right)\left( {1 - 2^{ - s - 1}} \right)\zeta \left( {s + 1}
\right) .\quad
\]

\vspace{3mm}

The term-by-term integration over any finite range is permissible
since the series $f\left( x \right)$ is boundedly convergent.

\vspace{3mm}

The series on the left hand side is

\[
\frac{\pi ^s}{s}\left[ {1 + \sum\limits_{m = 1}^\infty {\left( { -
1} \right)^m\left\{ {\left( {m + 1} \right)^s - m^s} \right\}} }
\right].
\]

\vspace{3mm}

This series is convergent for $s < 1$ and also uniformly convergent
for $\Re \left( s \right) < 1$.

\vspace{3mm}

Its sum is therefore an analytic function of $s$, regular for $\Re
\left( s \right) < 1$.

\vspace{3mm}

But for $s < 0$ it is $2\left( {1^s - 2^s + 3^s - ...} \right) =
2\left( {1 - 2^{s + 1}} \right)\zeta \left( { - s} \right).$

\vspace{3mm}

Its sum is therefore the same analytic function of $s$ for $\Re
\left( s \right) < 1$.

\vspace{3mm}

Hence, for $0 < s < 1$,

\[
\frac{\pi ^{s + 1}}{2s}\left( {1 - 2^{s + 1}} \right)\zeta \left( {
- s} \right) = \Gamma \left( s \right)\sin \left( {\frac{s\pi }{2}}
\right)\left( {1 - 2^{ - s - 1}} \right)\zeta \left( {s + 1}
\right),
\]

\vspace{3mm}

\noindent and the functional equation again follows.

\vspace{3mm}

We close this section, indicating that in the Appendix of this paper
we provide an application of Poisson's summation formula.

\vspace{6mm}

\begin{center}
3. THE SPECIAL FUNCTION $A\left( s \right)$\
\end{center}

\vspace{4mm}

At this stage, let us introduce a special function defined (in a
somewhat inelegant manner) by the following symbolic relation:

\begin{equation}
\label{eq3} A\left( {1 - s} \right) = \frac{\Gamma \left( {1 - s}
\right)\zeta \left( {1 - s} \right)L\left( {1 - s} \right)}{\pi ^{
{1 - s} }}.
\end{equation}

\vspace{1mm}

Changing (1) and (2) in (3), and considering Euler's complement
formula, that is true for the identification principle of the
relations among analytic functions, on the whole complex plane,
except for the integer values: $s = 0,\pm 1,\pm 2,\pm 3,\pm 4,\pm
5,\pm 6,\pm 7,....$, that is:

\[
\Gamma \left( s \right) \cdot \Gamma \left( {1 - s} \right) =
\frac{\pi }{\sin \left( {\pi s} \right)}
\]

\vspace{2mm}

\noindent we have:

\vspace{2mm}

\begin{equation}
\label{eq4} A\left( {1 - s} \right) = \frac{\Gamma \left( s
\right)\zeta \left( s \right)L\left( s \right)}{\pi ^{s }} = A\left(
s \right).
\end{equation}

\vspace{2mm}

The $A\left( s \right)$ function is actually a meromorphic function,
that satisfies the following remarkable identity:

\vspace{2mm}

\begin{equation}
\label{eq5} A\left( s \right) = A\left( {1 - s} \right).
\end{equation}

\vspace{2mm}

Such an identity shows a symmetry of the $A\left( s \right)$
function with respect to the vertical straight line  $\Re \left( s
\right) = 1 / 2$ [in particular let us consider the plot produced in
the interval $0 < \Re \left( s \right) < 1$ by the software product
\textit{DERIVE}\footnote{ \textit{DERIVE} is a Computer Algebra
System distribuited by Texas Instruments.}\textit{ Version 6.1 } for
Windows, for the real part of the analytic function $A\left( s
\right)$ ( see the curve in Fig. 1) ].

\medskip
\centerline{\epsfig{figure=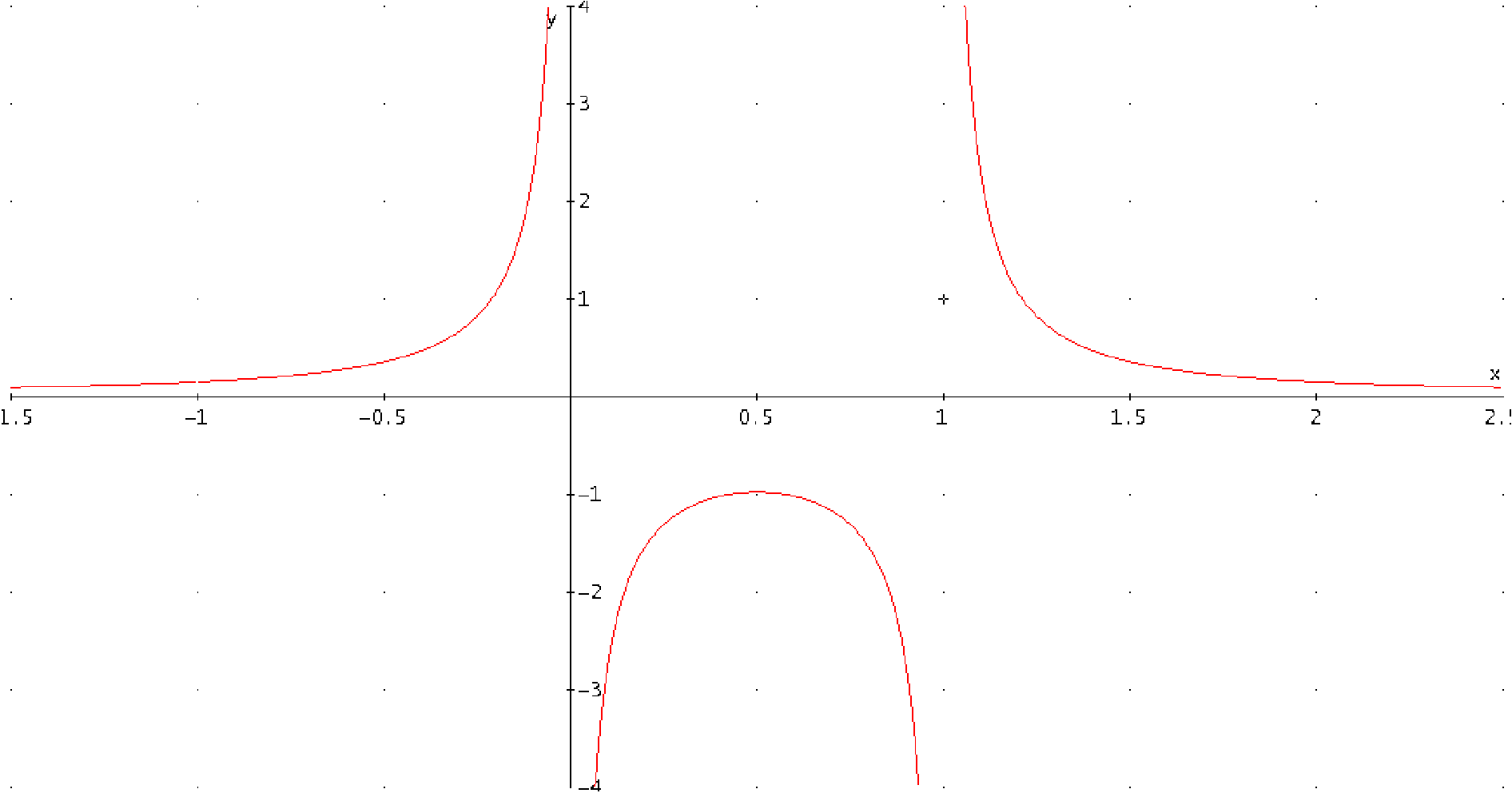,width=9cm,height=5.5cm}}
\medskip
\baselineskip=0.20in
\begin{center}
{\bf Fig. 1}
\end{center}

\vspace{6mm}

The identity (5) was obtained in a functional manner and leads to a
functional equation stemming from the entire function:

\begin{equation}
\label{eq6} \frac{\Gamma \left( s \right)  \zeta \left( s \right)
\,L\left( s \right)}{\Gamma \left( {1 - s} \right) \zeta \left( {1 -
s} \right)  \,L\left( {1 - s} \right)} = \;\frac{\pi ^s}{\pi ^{1 -
s}} = \exp \left[ {\left( {2s - 1} \right)\log \pi } \right].
\end{equation}

\vspace{3mm}

To verify this, we proceed as follows.

\vspace{3mm}

From (6) we get the ratio between $\zeta \left( s \right)$ and
$\zeta \left( {1 - s} \right)$:

\[
\frac{\zeta \left( s \right)}{\zeta \left( {1 - s} \right)\;} =
\;\frac{\pi ^s\Gamma \left( {1 - s} \right)L\left( {1 - s}
\right)}{\pi ^{1 - s}\Gamma \left( s \right)L\left( s \right)}.
\]

\vspace{3mm}

From this, re-using (2), we immediately get Riemann's well known
functional equation [4, pag. 35, formula (23)]:

\begin{equation}
\label{eq7} \zeta \left( s \right) = 2^s\,\pi ^{s - 1}\;\Gamma
\left( {1 - s} \right)\sin \left( {\frac{\pi s}{2}} \right)\;\zeta
\left( {1 - s} \right).
\end{equation}

\vspace{2mm}

In short, we have rewritten the following functional identity:

\begin{equation}
\label{eq8} \Gamma \left( s \right)\pi ^{ - s}\zeta \left( s
\right)L\left( s \right)\; = \;\Gamma \left( {1 - s} \right)\pi ^{ -
\left( {1 - s} \right)}\zeta \left( {1 - s} \right)L\left( {1 - s}
\right)\;
\end{equation}

\vspace{2mm}

\noindent in the form (7) and in this case we have not used Euler's
complement formula of the function $\Gamma \left( s \right)$, but we
have used only the  functional equation of $L\left( s \right)$.

\vspace{3mm}

In Riemann's memoir the functional identity is in the form\footnote{
$\prod {\left( s \right)} $ is related to the standard Gamma
function introduced by \textit{Legendre} by the equation $\prod
{\left( s \right)} = \Gamma \left( {s + 1} \right)$.\par }:

\[
\prod {\left( {\frac{s}{2} - 1} \right)} \;\pi ^{{ - s}
\mathord{\left/ {\vphantom {{ - s} 2}} \right.
\kern-\nulldelimiterspace} 2}\zeta \left( s \right) = \prod {\left(
{\frac{1 - s}{2} - 1} \right)} \;\pi ^{ - {\left( {1 - s} \right)}
\mathord{\left/ {\vphantom {{\left( {1 - s} \right)} 2}} \right.
\kern-\nulldelimiterspace} 2}\zeta \left( {1 - s} \right)
\]

\vspace{1mm}

\noindent that is slightly different from (\ref{eq8}).

\vspace{3mm}

This approach is the motivation for the title of paper : ``An
alternative form of the functional equation for Riemann's Zeta
function''.

\vspace{6mm}

\begin{center}
4. AN INTEGRAL REPRESENTATION
\end{center}

\vspace{4mm}

By using the identities [3, chap. X, p. 355, 10.15] :

\begin{equation}
\label{eq9} \Gamma \left( s \right)a^{ - s} = \int_0^\infty {x^{s -
1}e^{ - ax}dx \equiv M_s \left\{ {e^{ - ax}} \right\}} \
\end{equation}

\vspace{3mm}

\noindent where $M_s $ denotes the Mellin transform and

\begin{equation}
\label{eq10}\sum\limits_m {\left( { - 1} \right)} ^{m - 1}e^{ -
m^2x} = \frac{1}{2}\left[ {1 - \theta _4 \left( {0\left| {{ix}
\mathord{\left/ {\vphantom {{ix} \pi }} \right.
\kern-\nulldelimiterspace} \pi } \right.} \right)} \right]
 \
\end{equation}

\vspace{3mm}

\noindent where $\theta _4 \left( {z\left| \tau \right.} \right)$ is
the well-known $\theta \left( {z\left| \tau \right.} \right)$ theta
function of Jacobi [19, chap. XXI] and the summation variable $m$ is
to run over all positive integers, we derive the following integral
representation of the $A\left( s \right)$ function:

\begin{equation}
\label{eq11} A\left( s \right) = \frac{\pi ^{ - s}}{ {1 - 2^{1 - s}}
}\int_0^\infty {x^{s - 1}\left\{ {\frac{1}{4}\left[ {1 - \theta _4^2
\left( {0\left| {{ix} \mathord{\left/ {\vphantom {{ix} \pi }}
\right. \kern-\nulldelimiterspace} \pi } \right.} \right)} \right]}
\right\}dx}.
 \
\end{equation}

\vspace{4mm}

The key to obtain the previous formula are the works [6, pp.
409-410] and [7] by M.L. Glasser.

\vspace{3mm}

Combining the following two Mellin transforms:

\[\left( {1 - 2^{1 - 2s}} \right)\left[ \Gamma \left( s
\right)\zeta \left( {2s} \right) \right]
 = M_s \left\{
{\frac{1}{2}\left[ {1 - \theta _4 \left( {0\left| {{ix}
\mathord{\left/ {\vphantom {{ix} \pi }} \right.
\kern-\nulldelimiterspace} \pi } \right.} \right)} \right]}
\right\}\;,\quad \Re \left( s \right) > 0\]

 \noindent and \footnote{ See the correct formula in  \textit{"Solving some problems of advanced analytical nature
 posed in the SIAM-review"},
  by C.C. Grosjean, pag 432, Bull. Belg. Math. Soc. \textbf{3} (1996) and also in [7].}

\[\left( {1 - 2^{1 - 2s}} \right)\left[ {\Gamma \left( s
\right)\zeta \left( {2s} \right)} \right] - \left( {1 - 2^{1 - s}}
\right)\left[ {\Gamma \left( s \right)\zeta \left( s \right)L\left(
s \right)} \right] = M_s \left\{ {\frac{1}{4}\left[ {1 - \theta _4
\left( {0\left| {{ix} \mathord{\left/ {\vphantom {{ix} \pi }}
\right. \kern-\nulldelimiterspace} \pi } \right.} \right)}
\right]^2} \right\}\;,\quad \Re \left( s \right) > 0
 \]

\noindent the former is immediately obtained from Eqs. (9) and (10)
and the latter is obtained integrating term by term the following
remarkable identity, obtained from an identity\footnote{ See several
identities discovered by Jacobi [10].} by Jacobi and the result
$\theta _4^2 \left( {0\left| \tau \right.} \right) = {2k'K}
\mathord{\left/ {\vphantom {{2k'K} \pi }} \right.
\kern-\nulldelimiterspace} \pi $ ( [19], p. 479):

\[
\frac{1}{4}\left[ {1 - \theta _4^2 \left( {0\left| {{ix}
\mathord{\left/ {\vphantom {{ix} \pi }} \right.
\kern-\nulldelimiterspace} \pi } \right.} \right)} \right] =
\sum\limits_l {\left( { - 1} \right)} ^{{\left( {l - 1} \right)}
\mathord{\left/ {\vphantom {{\left( {l - 1} \right)} 2}} \right.
\kern-\nulldelimiterspace} 2}\left[ {e^{l x} + 1} \right]^{ - 1}
\]

\vspace{1mm}

\noindent(here the sum is to be expanded as a geometric series in
$e^{ - l x}$:
\[
e^{ - l x} - e^{ - 2l x} + e^{ - 3l x} - e^{ - 4l x} + e^{ - 5l x} -
\ldots = \left[ {e^{ l x} + 1} \right]^{ - 1}
\]

and the summation variable $l$ is to run over all positive odd
integers), we are in the position to determine the integral
representation (11) for the $A\left( s \right)$ function, by the
linearity property of the Mellin transformation, from the following
remarkable identity:

\[
A\left( s \right) = \frac{\Gamma \left( s \right)\zeta \left( s
\right)L\left( s \right)}{\pi ^s } \equiv \frac{\pi ^{ - s}}{ {1 -
2^{1 - s}} }\left[ M_s \left\{ {\frac{1}{2}\left[ {1 - \theta _4
\left( {0\left| {{ix} \mathord{\left/ {\vphantom {{ix} \pi }}
\right. \kern-\nulldelimiterspace} \pi } \right.} \right)} \right]}
\right\}\; -  M_s \left\{ {\frac{1}{4}\left[ {1 - \theta _4 \left(
{0\left| {{ix} \mathord{\left/ {\vphantom {{ix} \pi }} \right.
\kern-\nulldelimiterspace} \pi } \right.} \right)} \right]^2}
\right\} \right] \]

\[=\frac{\pi ^{ - s}}{ {1 - 2^{1 - s}} }M_s \left\{
{\frac{1}{4}\left[ {1 - \theta _4^2 \left( {0\left| {{ix}
\mathord{\left/ {\vphantom {{ix} \pi }} \right.
\kern-\nulldelimiterspace} \pi } \right.} \right)} \right]}
\right\}.
\]

\vspace{2mm}

The integral representation (11) is valid for $ \Re \left( s \right)
> 0$, in fact we can see that it contains the $\;\;$ product $\left( {1
- 2^{1 - s}} \right)\zeta \left( s \right)$ and it gives origin to
Euler's function $\eta \left( s \right)$ [known as Dirichlet's Eta
function], which is defined just for $\Re \left( s \right) > 0$
through the following alternating series:

\[\,
\eta \left( s \right) = \sum\limits_{n = 1}^\infty {\;\frac{\left( {
- 1} \right)^{n - 1}}{\left( n \right)^s}}.
\]
\vspace{8mm}

\begin{center}
5. THE $A\left( s \right)$\ FUNCTION AND THE HARMONIC SUM
\end{center}

\vspace{4mm}

Define
\begin{equation}
\label{eq12} H\left( x \right) = \sum\limits_{n = 1}^\infty
\sum\limits_{m = 0}^\infty  {\left( { - 1} \right)^m {e^{ - n\left(
{2m + 1} \right)x}} }. \
\end{equation}

\vspace{2mm}

\ This sum is also called the ``harmonic sum''.

\vspace{2mm}

Applying the basic functional properties of the Mellin transform (
[11], see APPENDIX), we find:

\begin{equation}
\label{eq13} \ M_s \left\{ {H\left( x \right)} \right\} =
\int\limits_0^\infty {H\left( x \right)} \;x^{s - 1}dx =
\sum\limits_{n = 1,m = 0}^\infty {\left( { - 1}
\right)^m\int\limits_0^\infty {e^{ - n\left( {2m + 1} \right)x}} }
\;x^{s - 1}dx =
\end{equation}

\[\sum\limits_{n = 1,m = 0}^\infty {\left( { - 1}
\right)^mn^{ - s}\left( {2m + 1} \right)^{ - s}\int\limits_0^\infty
{e^{ - t}} } \;t^{s - 1}dt = \sum\limits_{n = 1}^\infty {n^{ -
s}\sum\limits_{m = 0}^\infty {\left( { - 1} \right)^m} \left( {2m +
1} \right)^{ - s}\int\limits_0^\infty {e^{ - t}} } \;t^{s - 1}dt  =
\
\]

\vspace{2mm}

\[
\zeta \left( s \right)L\left( s \right)\Gamma \left( s \right).
\]

\vspace{5mm}

The interchange of summation and integration is legitimate by
Fubini's theorem.

\vspace{3mm}

Now we observe that:

\[
\sum\limits_{n = 1}^\infty \sum\limits_{m = 0}^\infty {\left( { - 1}
\right)^m {e^{ - n\left( {2m + 1} \right)x}} } = \sum\limits_{n =
1}^\infty {\;\frac{e^{ - nx}}{1 + e^{ - 2nx}}}
\]

\newpage

\noindent and besides that:

\[
H\left( x \right) = \sum\limits_{n = 1}^\infty {\;\frac{e^{ - nx}}{1
+ e^{ - 2nx}}} \le \sum\limits_{n = 1}^\infty {e^{ - nx}} = {e^{ -
x}} \mathord{\left/ {\vphantom {{e^{ - x}} {\left( {1 - e^{ - x}}
\right)}}} \right. \kern-\nulldelimiterspace} {\left( {1 - e^{ - x}}
\right)}.
\]

\vspace{4mm}

This result estabilishes that the fundamental strip of the Mellin
transform is $\left( {1, + \infty } \right)$ and that (13) exits for
any complex number $s > 1$.

\vspace{4mm}

From (13) we obtain another integral representation for the $A\left(
s \right)$ function, that is~:

\begin{equation}
\label{eq14} A\left( s \right) = \frac{1}{\pi
^s}\int\limits_0^\infty {H\left( x \right)} \;x^{s - 1}dx \quad
\quad \quad \Re \left( s \right) > 1
\end{equation}

\noindent where $H\left( x \right)$ is the function defined by the
infinite harmonic sum (12).

\vspace{3mm}

Now we use (14) to give an independent proof of (5) that does not
use (1) and (2).

\vspace{3mm}

We first consider the following transformation of $H\left( x
\right)$\textbf{, }that is:

\vspace{3mm}

\begin{equation}
\label{eq15} H\left( x \right) = \frac{\pi }{4x} - \frac{1}{4} +
\frac{\pi }{x}H\left( {\frac{\pi ^2}{x}} \right).
\end{equation}

\vspace{3mm}

This result is an immediate consequence of the Entry 11 in
\textit{Ramanujan's Notebooks} [2, pag. 258]:

\vspace{5mm}

\textit{Let }$\alpha ,\beta > 0$\textit{ with }$\alpha \beta = \pi
$\textit{, and let} $n$ \textit{be real with }$\left| n \right| <
\beta \mathord{\left/ {\vphantom {\beta 2}} \right.
\kern-\nulldelimiterspace} 2.$ \textit{Then}

\[
\alpha \left\{ {\frac{1}{4}\sec \left( {\alpha n} \right) +
\sum\limits_{k = 1}^\infty {\;\chi \left( k \right)\frac{\cos \left(
{\alpha \,nk} \right)}{e^{\alpha ^2k} - 1}} } \right\} = \beta
\left\{ {\frac{1}{4} +{\frac{1}{2}} \sum\limits_{k = 1}^\infty
{\;\frac{\cosh \left( {2\beta \,nk} \right)}{\cosh \left( {\beta
^2k} \right)}} } \right\}
\]

\textit{where}

\[
\chi \left( k \right) = \left\{ {\begin{array}{l}
 \,0\,\quad \;for\;k\;even \\
 \;1\,\quad for\;k\; \equiv \;1\;\bmod \;4 \\
 -1\,\,\,\,for\;k\; \equiv \;3\;\bmod \;4\; . \\
 \end{array}} \right.
\]

\vspace{3mm}

In the Appendix, by making use of the Poisson summation formula and
a remarkable Fourier cosine transform, we give a proof of the
infinite series identity, which in the author's view is simpler than
those given in [2] by Berndt.

\vspace{3mm}

For $n = 0$ this reads

\begin{equation}
\label{eq16} \alpha \left\{ {\frac{1}{4} + \sum\limits_{k =
1}^\infty {\;\chi \left( k \right)\frac{1}{e^{\alpha ^2k} - 1}} }
\right\} = \beta \left\{ {\frac{1}{4} +{\frac{1}{2}} \sum\limits_{k
= 1}^\infty {\;\frac{1}{\cosh \left( {\beta ^2k} \right)}} }
\right\}.
\end{equation}

\vspace{3mm}

Replacing $\cosh \left( x \right)$ by the exponential functions,
expanding the geometric series and rearranging the sums we obtain

\begin{equation}
\label{eq17} {\frac{1}{2}}\sum\limits_{k = 1}^\infty
{\;\frac{1}{\cosh \left( {\beta ^2k} \right)}} = \sum\limits_{m =
1}^\infty {\;\chi \left( m \right)\frac{1}{e^{\beta ^2m} - 1}}.
\end{equation}

\vspace{3mm}

Now considering the definition of the harmonic sum $H\left( x
\right)$\textbf{,} we have

\begin{equation}
\label{eq18} H\left( x \right) = \sum\limits_{m = 0}^\infty {\left(
{ - 1} \right)^m\sum\limits_{n = 1}^\infty  {e^{ - n\left( {2m + 1}
\right)x}} = \sum\limits_{m = 0}^\infty {\;\chi \left( m
\right)\frac{1}{e^{mx} - 1}}}.
\end{equation}

\vspace{3mm}

So that, taking into account the constraint $\alpha \beta = \pi , $
we find with (16), (17) and (18):

\[
\alpha \left\{ {\frac{1}{4} + H\left( {\alpha ^2} \right)} \right\}
= \beta \left\{ {\frac{1}{4} + H\left( {\beta ^2} \right)} \right\}.
\]

\vspace{3mm}

Finally, we substitute $\alpha = \sqrt x  , \; \beta = \pi
\mathord{\left/ {\vphantom {\pi {\sqrt x }}} \right.
\kern-\nulldelimiterspace} {\sqrt x }$ in the above relation and
obtain the desired transformation (15).

\vspace{3mm}

Plugging this back into our integral (14), we get

\[
A\left( s \right) = \frac{1}{\pi ^s}\int\limits_0^\infty {H\left( x
\right)} \;x^{s - 1}dx
\]

\[
 = \frac{1}{\pi ^s}\int\limits_\pi ^\infty {H\left( x \right)} \;x^{s - 1}dx
+ \frac{1}{\pi ^s}\int\limits_0^\pi {H\left( x \right)} \;x^{s -
1}dx
\]

\[
 = \frac{1}{\pi ^s}\int\limits_\pi ^\infty {H\left( x \right)} \;x^{s - 1}dx
+ \frac{1}{\pi ^s}\int\limits_0^\pi {\left\{ {\frac{\pi }{4x} -
\frac{1}{4} + \frac{\pi }{x}H\left( {\frac{\pi ^2}{x}} \right)}
\right\}} \;x^{s - 1}dx
\]

\[
 = \frac{1}{\pi ^s}\int\limits_\pi ^\infty {H\left( x \right)} \;x^{s - 1}dx
+ \frac{1}{4} \;\frac{1}{s\left( {s - 1} \right)} + \frac{1}{\pi
^s}\int\limits_\pi ^\infty {\frac{\pi ^{2s - 1}}{u^s}} \;H\left( u
\right)du\quad \quad \quad \left( {u = \frac{\pi ^2}{x}} \right)
\]

\[
 = + \frac{1}{4} \;\frac{1}{s\left( {s - 1} \right)} + \frac{1}{\pi
^s}\int\limits_\pi ^\infty {H\left( x \right)} \;x^{s - 1}dx +
\frac{1}{\pi ^{1 - s}}\int\limits_\pi ^\infty {H\left( x \right)}
\;x^{ - s}dx.
\]

\vspace{3mm}

Now the whole expression is symmetrical under $s \mapsto 1 - s$, the
integrals on the right hand side define a holomorphic function for
all $s \in
 \textbf{\textit{C, }}$ and so (5) follows.

\newpage

We summarize this results also in the following theorem:

\vspace{2mm}

\textbf{Theorem 1}: \textit{The Function}

\[
\zeta \left( s \right) = \frac{\pi ^s}{\Gamma \left( s
\right)\;L\left( s \right)}\left\{ {\frac{1}{4\;s\left( {s - 1}
\right)} + \int\limits_\pi ^\infty \;H\left( x \right)\left[ {\left(
{\frac{x}{\pi }} \right)^s + \left( {\frac{x}{\pi }} \right)^{1 -
s}} \right]\;d\,log\;x} \right\}
\]

\vspace{2mm}

\textit{is meromorphic with a simple pole at }$s = 1$ \textit{with
residue 1}.

\vspace{2mm}

Here is the computation of residue for $s = 1$ is [ observe that
$L\left( 1 \right) = arctg\left( 1 \right) = \pi \mathord{\left/
{\vphantom {\pi 4}} \right. \kern-\nulldelimiterspace} 4$\,\textbf{
- } a direct consequence of the infinite series expansion of the
arctangent (the Madhava-Gregory series)]:

\[
\lim \;_{s \to 1} \left( {s - 1} \right)\;\zeta \left( s \right) =
\frac{\pi }{\Gamma \left( 1 \right)\;L\left( 1 \right)} \cdot
\frac{1}{4} = 1.
\]

\vspace{6mm}

\begin{center}
6. CONCLUSION
\end{center}

\vspace{3mm}

Riemann gave two proofs of the functional equation (and functional
identity) in his ground-breaking paper [14], the former argument
essentially consists in proving the meromorphic continuation of the
$\zeta \left( s \right)$ function and uses contour integration,
while the latter, conceptually more difficult, using the $ $ $
\theta _3 \left( {0\left| {ix} \right.} \right)$ theta function,
requires the Mellin transformation.

\vspace{6mm}

In particular Riemann obtains the symmetrical formulation of the
functional identity by using two of the basic identities of the
Gamma function, that are: Euler's reflection (or complement) formula
and Legendre's duplication formula [13], which was discovered in
1809 and was surely unknown to Euler:

\[
\Gamma \left( s \right)\;\Gamma \left( {s + \raise0.7ex\hbox{$1$}
\!\mathord{\left/ {\vphantom {1
2}}\right.\kern-\nulldelimiterspace}\!\lower0.7ex\hbox{$2$}} \right)
= \frac{\sqrt s }{2^{2s - 1}}\Gamma \left( {2s} \right).
\]

\vspace{3mm}

Here we have introduced the complex $A\left( s \right)$\textbf{
}function and we have established another symmetrical formulation of
the functional equation for the Riemann Zeta function by using the
reflection formulae of the $\zeta \left( s \right)$, $L\left( s
\right)$ and $\Gamma\left( s \right)$ functions, all well-known by
Euler.

\vspace{3mm}

Using the definition of the $A\left( s \right)$ function we are also
able to obtain several integral representations of $A\left( s
\right)$ (with $\Re \left( s \right) > 0$ or $\Re \left( s \right)
> 1)$, that connect in an amazing way the $\zeta \left( s
\right)$ function with the independent transcendent $L\left( s
\right)$ function.

\vspace{2mm}

In addition, as a direct consequence of the second integral
representation, we have obtained an analytic continuation of the
same function by the Mellin transform of a function, defined by an
infinite harmonic sum, and using an identity of Ramanujan.

\vspace{2mm}

This last result represents another proof of the functional equation
for the $A\left( s \right)$ function, that is independent of the
three reflection formulae of the $\zeta \left( s \right)$, $L\left(
s \right)$ and $\Gamma\left( s \right)$ functions.

\vspace{2mm}

Finally, it is possible to state the following theorem:

\vspace{2mm}

\textbf{Theorem 2}: \textit{The} $A\left( s \right)$ \textit{
function extends itself as a meromorphic function in the complex
field }\textbf{\textit{C, }}\textit{in a regular way, except for the
simple poles at s = 0,1 [ respectively determined by the Gamma
function} $\Gamma \left( s \right)$\textit{ and by the Zeta function
}$\zeta \left( s \right)$ ] \textit{and satisfies the remarkable
functional equation}:

\[
A\left( s \right) = A\left( {1 - s} \right).
\]

\vspace{2mm}

\textit{As the complex zeros of the} $A\left( s \right)$
\textit{function coincide with the nontrivial zeros of the }$\zeta
\left( s \right)$\textit{ and} $L\left( s \right)$
\textit{functions, they are localized in the strip, determined by
}$0 \le \Re \left( s \right) \le 1$.

\vspace{2mm}

Let us keep in mind that in the $A\left( s \right)$ function the
singularities of the  $\Gamma \left( s \right)$ function, that we
can find in the negative real axis are cancelled by trivial zeros of
the two Euler's Zeta functions $ \zeta \left( s \right)$ and
$L\left( s \right)$.

\vspace{1mm}

It is in fact immediate to verify, from each of the functional
equations (1) and (2), exploiting the zeros of the trigonometric
functions \textit{cosine} and \textit{sine}, that:

\vspace{2mm}

 $\;\;\;\;\;\;\;\;\;\;\;\;\;\zeta \left( s \right) = 0\;\;$  for  $s = - 2, - 4, - 6,
- 8,...$ $ $ and $   $ $L\left( s \right) = 0\;\;$  for  $s = - 1, -
3, - 5, - 7,...$

\vspace{2mm}

Since the Gamma function has no zeros, and since the Dirichlet Beta
function and the Riemann Zeta function have an Euler product (see \S
1. Introduction), which shows that both are nonvanishing in the
right half plane $\Re \left( s \right) > 1$, the function $A\left( s
\right) = \Gamma \left( s \right){\kern 1pt} {\kern 1pt} \pi ^{ -
s}\zeta \left( s \right)L\left( s \right)$ has no zeros in $\Re
\left( s \right) > 1$.

\vspace{2mm}

By the functional equation $A\left( s \right) = A\left( {1 - s}
\right)$, it also has no zeros in $\Re \left( s \right) < 0$: thus
all the zeros have their real parts between 0 and 1.

\vspace{2mm}

Moreover, if all the complex zeros of the function $A\left( s
\right)$ have their real part equal to
$\raise.5ex\hbox{$\scriptstyle 1$}\kern-.1em/
\kern-.15em\lower.25ex\hbox{$\scriptstyle 2$} $, we shall obtain, as
results, both a proof of the Riemann Hypothesis, and the following
assertion [9] of Tschebyschef:

\vspace{2mm}

\quad\textit{       The function} \quad $F    \left( y \right) = e^{
- 3y} - e^{ - 5y} + e^{ - 7y} + e^{ - 11y} - .... = \sum\limits_{p
\succ 2} {\left( { - 1} \right)^{\frac{p + 1}{2}}} e^{ - py}$

\noindent \quad \textit{tends to infinity, as} $y \to 0$
$\quad\quad$ (the summation variable $p$ is to run over all odd
primes).

\vspace{2mm}

Indeed, in the paper [9] Hardy and Littlewood  prove that the
statement made by Tschebyschef is true if all complex zeros of the
function $L\left( s \right)$ have their real part equal to
$\raise.5ex\hbox{$\scriptstyle 1$}\kern-.1em/
\kern-.15em\lower.25ex\hbox{$\scriptstyle 2$} $.

\vspace{2mm}

This result confirms, in a very subtle way, the preponderance of
primes of the form $4m + 3$ [16, pag. 125].

\vspace{2mm}

Roughly speaking, there are "more" primes congruent to 3 mod 4 than
congruent to 1 mod 4.

\vspace{2mm}

The historical memoir of Riemann on the Zeta function has been
naturally extended to the family of Dirichlet $L$-functions
including the Riemann hypothesis.

\vspace{2mm}

The so-called Grand Riemann Hypothesis asserts that all the zeros of
$L\left( {s,\chi } \right)$ in the critical strip $0 < \Re \left( s
\right) < 1$ are on the critical line $\Re \left( s \right) =
\frac{1}{2}$ and this is a good reason to finish the paper right
here.

\vspace{2mm}

\begin{center}
7. APPENDIX
\end{center}

\vspace{2mm}

We prove the following infinite series identity of Ramanujan (
written in the inverse order):

\begin{equation}
\label{eq19} {\beta \left\{ {\frac{1}{4} + \sum\limits_{k =
1}^\infty {\;\frac{\cosh \left( {2\beta \,nk} \right)}{\cosh \left(
{\beta ^2k} \right)}} } \right\} = \alpha \left\{ {\frac{1}{4}\sec
\left( {\alpha n} \right) + \sum\limits_{k = 1}^\infty {\;\chi
\left( k \right)\frac{\cos \left( {\alpha \,nk} \right)}{e^{\alpha
^2k} - 1}} } \right\}}
\end{equation}

\vspace{3mm}

\noindent with $\alpha ,\beta > 0$\textbf{, }$\alpha \beta = \pi $
and $n$ real with $\left| n \right| < \beta \mathord{\left/
{\vphantom {\beta 2}} \right. \kern-\nulldelimiterspace} 2,$ using
the Poisson summation formula:

\begin{equation}
\label{eq20} {\frac{1}{2}f\left( 0 \right) + \sum\limits_{k =
1}^\infty {f\left( k \right)\; = \int\limits_0^\infty {f\left( x
\right)dx} } + 2\sum\limits_{k = 1}^\infty {\;\int\limits_0^\infty
{f\left( x \right)\cos \left( {2k\pi \,x} \right)} dx}}
\end{equation}

\vspace{2mm}

\noindent and the following Fourier cosine transform [3, cap. VII,
pag.174, 7.112] with $0 < a < b$:

\begin{equation}
\label{eq21} { F_c \left( u \right) = \;\int\limits_0^\infty
{\frac{\cosh \left( {at} \right)}{\cosh \left( {bt} \right)}\cos
\left( {ut} \right)} dt\; = \;\frac{\pi }{b}\frac{\cos \left(
{\frac{\pi \,a}{2b}} \right)\cosh \left( {\frac{\pi \,u}{2b}}
\right)}{\cos \left( {\frac{\pi \,a}{b}} \right) + \cosh \left(
{\frac{\pi \,u}{b}} \right)}}.
\end{equation}

\vspace{2mm}

We take into (20):\[ f\left( x \right) = \frac{\cosh \left( {2\beta
{\kern 1pt} n{\kern 1pt} x} \right)}{\cosh \left( {\beta ^2x}
\right)}
\]

\vspace{2mm}

\noindent and considering (21) with $a = 2\beta {\kern 1pt} n,\;b =
\beta ^2,\;u = 2k{\kern 1pt} \pi $ we have then:

\[
\;\int\limits_0^\infty {\frac{\cosh \left( {2\beta {\kern 1pt}
n{\kern 1pt} x} \right)}{\cosh \left( {\beta ^2x} \right)}\cos
\left( {2\pi {\kern 1pt} k{\kern 1pt} x} \right)} dx\; =
\;\frac{\alpha }{\beta }\frac{\cos \left( {n\alpha } \right)\cosh
\left( {k\alpha ^2} \right)}{\cos \left( {2n\alpha } \right) + \cosh
\left( {2k\alpha ^2} \right)}.
\]

\vspace{3mm}

By the relation $\cos z = \cosh \left(i z \right)$ and the following
prostapheresis-formula:

\[
\cosh \left( p \right) + \cosh \left( q \right) = \;2\cosh \left(
{\frac{p + q}{2}} \right) \cdot \cosh \left( {\frac{p - q}{2}}
\right)
\]

\vspace{1mm}

\noindent we have immediately:

\begin{equation}
\label{eq22} { \;\int\limits_0^\infty {\frac{\cosh \left( {2\beta
{\kern 1pt} n{\kern 1pt} x} \right)}{\cosh \left( {\beta ^2x}
\right)}\cos \left( {2\pi {\kern 1pt} k{\kern 1pt} x} \right)} dx\;
= \;\frac{\alpha }{4\beta }\left[ {sech \left( {k\alpha ^2 -
in\alpha } \right) + sech \left( {k\alpha ^2 + in\alpha } \right)}
\right]}.
\end{equation}

\newpage

Observe that setting $k = 0$ we have:

\begin{equation}
\label{eq23} {\;\int\limits_0^\infty {\frac{\cosh \left( {2\beta
{\kern 1pt} n{\kern 1pt} x} \right)}{\cosh \left( {\beta ^2x}
\right)}} dx\; = \;\frac{\alpha }{2\beta }\sec \left( {n\alpha }
\right)}.
\end{equation}

\vspace{1mm}

In short considering:

\[
\frac{1}{\cosh {\kern 1pt} \,\left( x \right)}\; = sech\,\left( x
\right)\; = 2\sum\limits_{r = 1}^\infty {\chi \left( r \right)}
\,e^{ - rx}
\]

\vspace{1mm}

\noindent and besides:

\[
\left[ {sech\left( {k\alpha ^2 - in\alpha } \right) + sech\left(
{k\alpha ^2 + in\alpha } \right)} \right] = \;4\sum\limits_{r =
1}^\infty {\chi \left( r \right)} \,e^{ - 2k\alpha ^2}\cos \left(
{\alpha \,n\,r} \right)
\]

\vspace{1mm}

\noindent we also deduce that:

\[
\sum\limits_{k = 1}^\infty {\left[ {sech\left( {k\alpha ^2 -
in\alpha } \right) + sech\left( {k\alpha ^2 + in\alpha } \right)}
\right] = } \;4\sum\limits_{r = 1}^\infty {\chi \left( r \right)}
\,\frac{\cos \left( {\alpha \,n\,r} \right)}{e^{\alpha ^2\,r} - 1}.
\]

\vspace{1mm}

As a consequence of (22) putting the above equality with (23) into
(20) we have the identity (19).

\vspace{3mm}

\begin{center}
8. ADDITIONAL REMARKS
\end{center}

\vspace{3mm}

\textit{Remark} 1. The functional equation (5) for $A\left( s
\right)$ is not new.

\vspace{3mm}

In [L. Lorenz, {\it Tidskr. Mat.} \textbf{1}, 97 (1871)] it is shown
that:

\vspace{2mm}

\[A\left( s \right) = \frac{\Gamma \left( s \right)}{4\pi
^s}\,\;Z\left| {\mathop {\begin{array}{l}
 0 \\
 \mathop 0\limits^{ \to } \\
 \end{array}}\limits^ \to } \right|\left( {1;\;2s} \right)
\]

\noindent in terms of Epstein's Zeta-function:

\[ Z\left( {1;2s}
\right) = \sum\limits_{\begin{array}{l}\;\;\;m,n = - \infty \\
 \left( {m,n} \right) \ne \left( {0,0} \right) \\
 \end{array}}^\infty \frac{1}{\left( {m^2 + n^2} \right)^s}
\]

\vspace{1mm}

[P. Epstein, {\it Zur Theorie allgemeiner Zetafunktionen. I.}, Math.
Ann. \textbf{56}, 615 (1903)].

\vspace{3mm}

This function satisfies the functional equation:
\[
\pi ^{-s}\Gamma \left( s \right)Z\left( {1;\;2s} \right)\; = \;\pi
^{ - \left( {1 - s} \right)}\Gamma \left( {1 - s} \right)Z\left(
{1;\;2 - 2s} \right)
\]
which is the same as $ A\left( s \right) = A\left( {1 - s} \right)$.

\vspace{3mm}

\textit{Remark} 2. One motivation for the surprisingly quick proof
of the symmetrical formulation [see (\ref{eq5}) or (\ref{eq8})] of
the functional equation of the Zeta function is that Euler himself
could have proved this remarkable identity with the three reflection
formulae (this is the reason for a dedication to Leonhard Euler).

\vspace{3mm}

\textit{Remark} 3. We could have obtained an analytic continuation
of the function $A\left( s \right)$ also from the first integral
representation (11), following a similar method (the transformation
law of theta function ) to the one used by Riemann, but we have
chosen a second opportunity, just to give a more innovative proof
with the identity of Ramanujan and therefore slightly different from
the classical one, that we find in Riemann's original memoir.


\vspace{2mm}

\begin{thebibliography}{99}

\vspace{4mm}


\bibitem{1} N.H. Abel, {\it Solution de quelques probl\`{e}mes \`{a}
l'aide d'integrales d\'{e}finies}, Mag. Naturvidenskaberne
\textbf{2}, 1823.

\vspace{2mm}

\bibitem{2} B. C. Berndt, {\it Ramanujan's Notebooks. Part II}, Springer-Verlag, New York, 1989.

\vspace{2mm}

\bibitem{3} V. Ditkinev et A. Proudnikov, {\it Transformations Integrales e Calcul Op\`{e}rationnel
}, Ed. Mir, Mosca, 1982.

\vspace{2mm}

\bibitem{4} I. Erdelyi et al. (ed), {\it Higher Trascendental Functions}, vol 1,
(Bateman Manuscript Project), McGraw-Hill Book Company, Inc., New
York, 1953.

\vspace{2mm}

\bibitem{5} S.R. Finch, {\it Mathematical Constants}, Cambridge Univ. Press, Cambridge, 2003.

\vspace{2mm}

\bibitem{6} M. L. Glasser, {\it The evaluation of lattice sums. I.: Analytic procedures},
        J. Math. Phis., Vol. \textbf{14}, March 1973.

\vspace{2mm}


\bibitem{7} M. L. Glasser and I.J. Zucker, in  {\it Theoretical Chemistry: Advances and Perspectives},
        Vol. \textbf{5}, D. Henderson and H. Eyring, Acad. Press, New York, pp. 67-139, 1980.

\vspace{2mm}

\bibitem{8} G. H. Hardy, {\it Divergentes Series}, Cambridge, 1949.

\vspace{2mm}


\bibitem{9} G. H. Hardy and J. E. Littlewood, {\it Contributions to the theory of the Riemann zeta-function and the theory of the distribution of primes},
Acta Mathematica, \textbf{41}, pp. 119-198, 1918.

 \vspace{2mm}

\bibitem{10} C. G. I. Jacobi, {\it Fundamenta Nova Theoriae Functionum Ellipticarum},(K\"{o}nigsberg, 1829), Sec. \textbf{40}.

\vspace{2mm}

\bibitem{11} P. Kirschenhofer, H. Prodinger, and W. Szpankowski, {\it Multidimensional digital searching and some new parameters in tries}, Tecnical Report
CSD-TR-91-052, Purdue University, July 1991.

\vspace{2mm}

\bibitem{12}M. Knopp and S. Robins, {\it Easy proofs of Riemann's functional equation for }$\zeta \left( s \right)$\textit{ and of Lipschitz summation},
Proc. Amer. Math. Soc. Vol. \textbf{129, }no.\textbf{ 7}, pp.
1915-1922, 2001.

\vspace{2mm}

\bibitem{13} A. M. Legendre, {\it M\'{e}moires de la classe des sciences
math\'{e}matiques et phisiques de l'Institut de France}, Paris, pp.
477-490, 1809.

\vspace{2mm}

\bibitem{14} B. Riemann, {\it Ueber die Anzahl der Primzahlen unter einer gegebenen  Gr\"{o}sse}, 1859.
In {\it Gesammelte Werke}, Teubner, Leipzig, 1892, Reprinted Dover,
New York, 1953.

\vspace{2mm}

\bibitem{15} H.M. Srivastava, J. Choi, {\it Series Associated with the Zeta and Related Functions}, Kluwer Academic Publishers,
Dordrecht, Boston, and London, 2001.

\vspace{2mm}

\bibitem{16} G. Tenenbaum et M. Mend\`{e}s France, {\it Les Nombres Premiers}, Presses Universitaires de France, Paris, 1997.

\vspace{2mm}

\bibitem{17} E. C. Titchmarsh and D. R. Heath-Brown, {\it The Theory of the Riemann Zeta-Function}, 2nd ed. Oxford,
England: Oxford University Press, 1986.

\vspace{2mm}

\bibitem{18} V. S. Varadarajan, {\it Euler and his work of infinite series,} Bulletin of the American Mathematical Society, vol.
\textbf{44}, no. \textbf{4}, pp. 515-539, 2007.

\vspace{2mm}

\bibitem{19} E.T. Whittaker, G.N. Watson, {\it A course of modern analysis}, 4th ed.,
         Cambridge University Press, Cambridge, 1988.


\end{thebibliography}
\end{document}